\documentclass{article}
\usepackage{amsmath}[1996/11/01]
\usepackage{amssymb,amsthm,amsxtra}
\usepackage{epic,eepic}
\usepackage{epsfig}

\def\[#1\]{\begin{equation}#1\end{equation}}
\makeatletter
\def\beq{%
   \relax\ifmmode
      \@badmath
   \else
      \ifvmode
         \nointerlineskip
         \makebox[.6\linewidth]%
      \fi
      $$
   \fi
}
\def\eeq{%
   \relax\ifmmode
      \ifinner
         \@badmath
      \else
         $$
      \fi
   \else
      \@badmath
   \fi
   \ignorespaces
}

\def\enddisplaymath{\eeq\global\@ignoretrue}
\makeatother

\newtheorem{thm}{Theorem}

\newtheorem{lem}[thm]{Lemma}
\newtheorem{prop}[thm]{Proposition}
\theoremstyle{remark}
\newtheorem*{rem}{Remark}
\newtheorem{rems}{Remark}[thm]

\theoremstyle{definition}
\newtheorem{defn}{Definition}

\numberwithin{equation}{section}
\numberwithin{thm}{section}

\DeclareMathOperator{\Tr}{Tr}
\DeclareMathOperator{\Exp}{\mathbb E}

\renewcommand{\Re}{\operatorname{Re}}

\DeclareMathOperator{\Prob}{\mathbb P}

\DeclareMathOperator{\erf}{erf}

\newcommand{\C}{\mathbb C}
\newcommand{\R}{\mathbb R}
\newcommand{\Z}{\mathbb Z}
\newcommand{\N}{\mathbb N}

\DeclareMathOperator{\FGOE}{F_{GOE}}
\DeclareMathOperator{\FGUE}{F_{GUE}}
\DeclareMathOperator{\FGSE}{F_{GSE}}

\DeclareMathOperator{\F0}{F_0}

\newcommand\psymmU{%
\begin{picture}(1,1)(0,0)%
\allinethickness{0.5pt}%
\path(0,0)(0,1)(1,1)(1,0)(0,0)%
\end{picture}}
\newcommand\psymmUU{%
\begin{picture}(1,1)(0,0)%
\allinethickness{0.5pt}%
\path(0,0)(0,1)(1,1)(1,0)(0,0)%
\put(0.5,0.5){\makebox(0,0){$\cdot$}}%
\end{picture}}
\newcommand\psymmO{%
\begin{picture}(1,1)(0,0)%
\allinethickness{0.5pt}%
\path(0,0)(0,1)(1,1)(1,0)(0,0)%
\path(0,0)(1,1)%
\end{picture}}
\newcommand\psymmS{%
\begin{picture}(1,1)(0,0)%
\allinethickness{0.5pt}%
\path(0,0)(0,1)(1,1)(1,0)(0,0)%
\path(1,0)(0,1)%
\end{picture}}
\newcommand\psymmu{%
\begin{picture}(1,1)(0,0)%
\allinethickness{0.5pt}%
\path(0,0)(0,1)(1,1)(1,0)(0,0)%
\path(0,0)(1,1)%
\path(0,1)(1,0)%
\end{picture}}

\newbox\tsymmUbox
\newbox\tsymmUUbox
\newbox\tsymmObox
\newbox\tsymmSbox
\newbox\tsymmubox
\setbox\tsymmUbox =\hbox{\kern0.75pt\setlength{\unitlength}{6pt}\psymmU \kern0.75pt}
\setbox\tsymmUUbox=\hbox{\kern0.75pt\setlength{\unitlength}{6pt}\psymmUU\kern0.75pt}
\setbox\tsymmObox =\hbox{\kern0.75pt\setlength{\unitlength}{6pt}\psymmO \kern0.75pt}
\setbox\tsymmSbox =\hbox{\kern0.75pt\setlength{\unitlength}{6pt}\psymmS \kern0.75pt}
\setbox\tsymmubox =\hbox{\kern0.75pt\setlength{\unitlength}{6pt}\psymmu \kern0.75pt}
\def\tsymmU{{\copy\tsymmUbox}}

\def\tsymmu{{\copy\tsymmubox}}
\newbox\symmUbox
\newbox\symmUUbox
\newbox\symmObox
\newbox\symmSbox
\newbox\symmubox
\setbox\symmUbox =\hbox{\kern0.75pt\setlength{\unitlength}{4.5pt}\psymmU \kern0.75pt}
\setbox\symmUUbox=\hbox{\kern0.75pt\setlength{\unitlength}{4.5pt}\psymmUU\kern0.75pt}
\setbox\symmObox =\hbox{\kern0.75pt\setlength{\unitlength}{4.5pt}\psymmO \kern0.75pt}
\setbox\symmSbox =\hbox{\kern0.75pt\setlength{\unitlength}{4.5pt}\psymmS \kern0.75pt}
\setbox\symmubox =\hbox{\kern0.75pt\setlength{\unitlength}{4.5pt}\psymmu \kern0.75pt}

\def\symmu{{\copy\symmubox}}

\begin{document}

\title{{\bf Limiting distributions for a polynuclear growth model \\
with external sources}} 
\author{{\bf Jinho Baik}\footnote{
Deparment of Mathematics, 
Princeton University, Princeton, New Jersey, 08544, 
jbaik@math.princeton.edu}
\footnote{Institute for Advanced Study, 
Princeton, New Jersey 08540}
 \ \ 
and \ \ {\bf Eric M. Rains}\footnote{AT\&T Research, New Jersey, 
Florham Park, New Jersey 07932, rains@research.att.com}}

\date{March 22, 2000}
\maketitle

\begin{abstract}

The purpose of this paper is to investigate the limiting distribution
functions for a polynuclear growth model with two external sources, which
was considered by Pr\"ahofer and Spohn in \cite{SpohnP1}.  Depending on the
strength of the sources, the limiting distribution functions are either the
Tracy-Widom functions of random matrix theory, or a new explicit function
which has the special property that its mean is zero.  Moreover, we obtain
transition functions between pairs of the above distribution functions in
suitably scaled limits.  There are also similar results for a discrete
totally asymmetric exclusion process.

\medskip


\end{abstract}


\section{Introduction}\label{sec-intro}

\begin{figure}[ht]
 \centerline{\epsfig{file=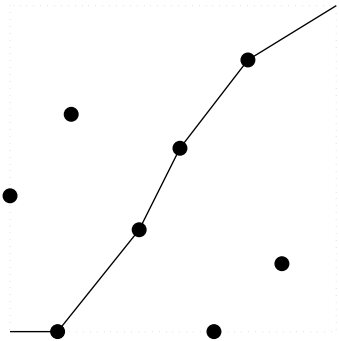, width=4.5cm}}
 \caption{A points selection and the longest up/right path}
\label{fig-example}
\end{figure}

Our main object of study is the following combinatorial problem.  Fix three
real parameters $t>0$ and $\alpha_\pm\ge 0$.  We construct a random set of
points in the unit square $[0,1]\times [0,1]\subset\R^2$, as follows.  Let
$P(\lambda)$ denote a Poisson variable of density $\lambda$.  We select
$P(t^2)$ points at random inside the square $(0,1)\times (0,1)$,
$P(\alpha_+t)$ points at random on the open bottom edge $(0,1)\times
\{0\}$, and $P(\alpha_-t)$ points at random on the open left edge
$\{0\}\times (0,1)$.  Hence no point is selected from the lower left vertex
and from the closed top and right edges.  For example, in Figure
\ref{fig-example}, 5 points are selected inside the square, and 2 and 1
points in the bottom and the left edges, respectively.  A (weakly) up/right
path is given by a sequence of points such that each point is (weakly)
above and to the right of its predecessor; thus the solid line in the
figure is a (weakly) up/right path from the lower left vertex to the upper
right vertex.
The length of a (weakly) up/right path is defined by the number of points
in the path; thus the up/right path of the example has length $4$.

We define $L(t)$ be the length of the longest (weakly) up/right path in
this random configuration of points.  (This, of course, implicitly depends
on the values of $\alpha_\pm$.)  E.g., the solid line in Figure
\ref{fig-example} is the longest up/right path of the example.  In general,
there may be more than one longest path, but we are only interested in the
length $L(t)$, hence the degeneracy is not an issue.  We note that the
above process can be thought of in an alternative way as a Poisson process
of intensity $1$ in the open $\R_+^2$, together with a Poisson process of
intensity $\alpha_+$ on the open half-line $\R_+\times\{0\}$ and a Poisson
process of intensity $\alpha_-$ on the open half-line $\{0\}\times \R_-$.
Then $L(t)$ is equal to the length of the longest (weakly) up/right path
from $(0,0)$ to $(t,t)$.

The main interest in this paper is the statistics of $L(t)$ as
$t\to\infty$, as this problem arose in a polynuclear growth (PNG)
model considered by Pr\"ahofer and Spohn in \cite{SpohnP1}.  PNG is a
simplified model for layer by layer growth in one spatial dimension.  At
each random nucleation position, an island of height 1 is formed and
spreads laterally with speed 1.  When two islands meet, they form one
island and keep spreading with the same speed.  The question is the
fluctuations of the height in the large time limit.

As a special case, suppose that a single island 
starts spreading at the origin and that further nucleation takes place 
only on top of this ground layer.
Furthermore, suppose that there are external nucleation sources 
at the two ends of the ground layer.
Then the height, for example, at the origin at time $\sqrt{2}t$ is 
equal to $L(t)$ in the above point selection process.
We refer the interested readers to the papers \cite{SpohnP1, SpohnP2} 
for details of the mapping from PNG to $L(t)$ and 
many other related works.


\medskip
In the above point selection process, a special case is when $\alpha_\pm=0$
: $P(t^2)$ points are selected inside the square with no points on the
edges.  (This corresponds to a PNG model without external nucleation
sources.)  It is interesting that in this case, there is a combinatorial
interpretation.  Consider a random permutation of $S_N$.  If we take $N$ as
Poisson variable $P(t^2)$, the longest increasing subsequence of a random
permutation of $S_N$ has the same statistics as $L(t)$ (see \cite{BDJ}
and references therein).  The limiting fluctuation of $L(t)$ in this case
is obtained in \cite{BDJ} : there is a distribution function $\FGUE(x)$
such that
\begin{equation}\label{eq1-1}
   \lim_{t\to\infty}\Prob\biggl(\frac{L(t)-2t}{t^{1/3}}\le x\biggr) 
= \FGUE(x), \qquad \alpha_\pm= 0.
\end{equation}
The convergences of all the moments are also proved in the same paper.
Here $\FGUE$ is the so-called 
\emph{GUE Tracy-Widom distribution function} 
of random matrix theory.  See Section \ref{sec:dist} for the explicit
formula of $\FGUE$, equation \eqref{eq2-8}.

The ``GUE'' in $\FGUE$ refers to the Gaussian unitary ensemble,
the set of all $N\times N$ Hermitian matrices together 
with the probability measure defined by 
\begin{equation}
  \frac1{Z_N}e^{-tr M^2}dM
\end{equation}
where $Z_N$ is a normalization constant
(see e.g. \cite{Mehta, deift}).
Tracy and Widom in \cite{TW1} 
proved that as $N\to\infty$, the properly centered and scaled 
largest eigenvalue of a random GUE matrix converges in distribution 
to $\FGUE$.
Therefore the above result \eqref{eq1-1} implies that 
in the limit, $L(t)$ with $\alpha_\pm=0$ 
and the largest eigenvalue of a random GUE matrix 
have the same fluctuations.

\medskip
In this paper, we obtain similar limiting distributions in the presence of
points on the edges.  When $\alpha_\pm>0$, the longest path begins by
following one of the edges for a while, then enters the square.  Hence when
$\alpha_\pm$ are small, the effect coming from edges is small, and we
expect $L(t)$ to have the same statistics as when $\alpha_\pm=0$, i.e. GUE
fluctuation.  But when $\alpha_\pm$ is large, the longest path lies mostly
on one of the edges, hence we expect Gaussian fluctuation.  It turns out
that the critical case is when $\alpha_\pm=1$.  We need to distinguish four
different cases to state the results.  In each case, $\alpha_\pm$ are
fixed.

\begin{itemize}
\item When $\alpha_\pm<1$, we obtain $\FGUE$ in the limit.
\item When either of $\alpha_\pm$ is greater than $1$, 
we obtain Gaussian fluctuation.
\item 
When one of $\alpha_\pm$ is equal to $1$ and the other is strictly less 
than $1$, we obtain 
$\FGOE(x)^2$.
In the above definition of the GUE, if we replace Hermitian matrices 
by real symmetric matrices, we obtain the
Gaussian orthogonal ensemble (GOE).
The limiting distribution of the (properly centered and scaled) 
largest eigenvalue of a random GOE matrix is 
given by $\FGOE$ in \cite{TW2} (see Section \ref{sec:dist} 
for the explicit formula). 
The above limit $\FGOE(x)^2$ can be interpreted as follows : 
take two random GOE matrices and superimpose their eigenvalues.
The largest of the superimposition of eigenvalues has limiting 
fluctuation $\FGOE^2$.
\item When $\alpha_\pm=1$, 
we have a new limiting distribution which we denote $\F0$,
for which we do not yet have a random matrix interpretation.
See Section \ref{sec:dist} for the definition and discussions.
\end{itemize}
More explicitly, we have the following theorem.

\begin{thm}\label{thm1}
For each fixed $\alpha_\pm$, as $t\to\infty$, we have the 
following results.
\begin{enumerate}
\item When $0\le \alpha_\pm \le 1$, 
\begin{equation}
  \lim_{t\to\infty}\Prob\biggr( \frac{L(t)-2t}{t^{1/3}}\le x \biggr) =
\begin{cases}
  \FGUE(x), \quad &0\le \alpha_\pm <1,\\
  \FGOE(x)^2, \quad &\text{$\alpha_+=1, 0\le \alpha_-<1$, or $\alpha_-=1, 
0\le\alpha_+<1$,}\\
  \F0(x), \quad &\alpha_\pm=1.
\end{cases}
\end{equation}
\item When at least one of $\alpha_\pm$ is greater than $1$, 
setting $\alpha=\max\{\alpha_+,\alpha_-\}$, 
\begin{equation}\label{e-i4}
  \lim_{t\to\infty} \Prob\biggl(\frac{L(t)-(\alpha+\alpha^{-1})t}
{\sqrt{\alpha-\alpha^{-1}}t^{1/2}}\le
x\biggr) =
\begin{cases}
  \erf(x), \quad &\alpha_+\neq \alpha_- \\
  \erf(x)^2, \quad &\alpha_+=\alpha_-.
\end{cases}
\end{equation}
\end{enumerate}
\end{thm}

\begin{rem}
As $\alpha$ tends to $\infty$, $L(t)$ tends to a variable with mean and
variance $\alpha t$; i.e., to $P(\alpha t)$, corresponding to the number of
points on the interval $(0,1)$.  For finite $\alpha$, the longest path
leaves the interval near $1-\alpha^{-2}$, then includes roughly
$2\alpha^{-1}t$ points inside the square, thus giving $L(t)\sim
P((\alpha-\alpha^{-1})t)+2\alpha^{-1}t$, agreeing with \eqref{e-i4}.
\end{rem}

\bigskip
The functions, $\FGUE$, $\FGOE$
and $\F0$, in the above theorem
are defined in Section \ref{sec:dist}, and we discuss their properties 
in the same section.

The above theorem shows that 
there are certain transitions around the points $\alpha_\pm=1$.
It is of interest to investigate these transitions in detail.
In Section \ref{sec:scale}, 
under proper scalings of $\alpha_\pm\to 1$, 
we obtain new classes of distribution 
functions interpolating the functions of the above theorem. 

In addition to the above point selection process, there is a closely
related lattice directed polymer model in 2-dimensional space, which is a
generalization of the model considered by Johansson in \cite{kurtj:shape}.
This model and a related exclusion process are discussed in Section
\ref{sec:excl}.

Finally 
the proofs of the theorems are discussed in Section \ref{sec:proof}.

\medskip
\noindent {\bf Notations.}
In many papers, 
the functions $\FGUE$ and $\FGOE$ are denoted by 
$F_2$ and $F_1$, respectively.

\medskip
\noindent {\bf Acknowledgments.}
The authors greatly appreciate many enlightening conversations and
communications with Michael Pr\"ahofer and Herbert Spohn 
who brought our interest to this problem.

\section{Limiting distribution functions}\label{sec:dist}

Let $u(x)$ be the solution to the Painlev\'e II (PII) equation,
\begin{equation}\label{as10}
  u_{xx}=2u^3+xu,
\end{equation}
with the boundary condition
\begin{equation}\label{as11}
  u(x)\sim -Ai(x)\quad\text{as}\quad x\to +\infty,
\end{equation}
where $Ai$ is the Airy function.
The proofs of the (global) existence and the uniqueness of the solution
were established in \cite{HM}.
The asymptotics as $x\to -\infty$ are given by 
(see e.g. \cite{HM,DZ2})
\begin{eqnarray}
\label{as17}  u(x) &=&
-Ai(x) + O\biggl( \frac{e^{-(4/3)x^{3/2}}}{x^{1/4}}\biggr),
\qquad\text{as $x\to +\infty$,}\\
\label{as18}  u(x) &=&
-\sqrt{\frac{-x}{2}}\biggl( 1+O\bigl(\frac1{x^2}\bigr)\biggr),
\qquad\text{as $x\to -\infty$.}
\end{eqnarray}
Recall that $Ai(x)\sim \frac{e^{-(2/3)x^{3/2}}}{2\sqrt{\pi}x^{1/4}}$
as $x\to +\infty$.
Define
\begin{equation}\label{as12}
  v(x):= \int_{\infty}^{x} (u(s))^2 ds,
\end{equation}
so that $v'(x)= (u(x))^2$.
We note another expression 
\begin{equation} 
  v(x)=u(x)^4+xu(x)^2-(u'(x))^2
\end{equation} 
which can be obtained by noting that the difference 
(i) has derivatives equal to zero by the PII equation, 
and (ii) becomes zero as $x\to\infty$ by \eqref{as17}.

The Tracy-Widom distribution functions are defined in terms of $u$ and $v$.
\begin{defn}[TW distribution functions]\label{def1}
Set
\begin{eqnarray}
\label{as46-1}   F(x) &:=& \exp\biggl(\frac12\int_x^{\infty} v(s)ds\biggr)
= \exp\biggl(-\frac12\int_x^{\infty} (s-x)(u(s))^2ds\biggr),\\
\label{as47-1}   E(x) &:=& \exp\biggl(\frac12\int_x^{\infty} u(s)ds\biggr),
\end{eqnarray}
and set
\begin{eqnarray}
\label{eq2-8}
   \FGUE(x) :=& F(x)^2 &=\exp\biggl(-\int_x^{\infty} (s-x)(u(s))^2ds\biggr),\\
\label{eq2-9}
   \FGOE(x) :=& F(x)E(x) &= \bigl(F_2(x)\bigr)^{1/2}
e^{\frac12\int_x^{\infty} u(s)ds} ,\\
\label{eq2-10}
   \FGSE(x) :=& F(x)\bigl\{E(x)^{-1}+E(x)\bigr\}/2
&= \bigl(F_2(x)\bigr)^{1/2}
\biggl[ e^{-\frac12\int_x^{\infty} u(s)ds} +
e^{\frac12\int_x^{\infty} u(s)ds} \biggr] /2 .
\end{eqnarray}
\end{defn}

It is proved by Tracy and Widom in \cite{TW1, TW2} 
that under proper centering
and scaling, the distribution of the largest eigenvalue of
a random GUE/GOE/GSE matrix converges to $\FGUE(x)$ / $\FGOE(x)$ / $\FGSE(x)$
as the size of the matrix becomes large.
The readers are referred to \cite{Mehta, deift} for definitions of various 
random matrix ensembles and their basic properties.
We note that from the asymptotics \eqref{as17} and \eqref{as18},
for some positive constant $c$,
\begin{eqnarray}
\label{as46}   F(x) &=& 1+ O\bigl( e^{-cx^{3/2}}\bigr)
\qquad \text{as $x\to+\infty$,}\\
\label{as47}   E(x) &=& 1+ O\bigl( e^{-cx^{3/2}}\bigr)
\qquad \text{as $x\to+\infty$,}\\
\label{as48}   F(x) &=& O\bigl( e^{-c|x|^{3}}\bigr)
\qquad\qquad \text{as $x\to-\infty$,}\\
\label{as49}   E(x) &=& O\bigl( e^{-c|x|^{3/2}}\bigr)
\quad\qquad \text{as $x\to-\infty$.}
\end{eqnarray}
Hence in particular, as $x\to+\infty$, all the above three functions become
$1$, and as $x\to -\infty$, they become $0$.  Monotonicity follows from the
fact that they are limits of sequences of distribution functions, and
therefore \eqref{eq2-8}-\eqref{eq2-10} are indeed distribution functions.

We need a new distribution function for the case when $\alpha_\pm=1$ in 
Theorem \ref{thm1}.
\begin{defn}
Set 
\begin{equation}\label{eq-f0}
  \F0(x)= \bigl\{1-(x+2u'(x)+2u(x)^2)v(x)\bigr\} \bigl(E(x)\bigr)^4 \FGUE(x).
\end{equation}
\end{defn}
The asymptotics \eqref{as17}-\eqref{as12} imply that $\F0(x)$ has limit 
$1$ as $x\to +\infty$ and $0$ as $x\to -\infty$.
The monotonicity of $\F0$ follows from the fact that 
it is a limit of distribution functions in Theorem \ref{thm1}.
It would of interest to have random matrix interpretation of the function
$\F0$, but so far we have been unable to identity $\F0$ as a quantity
arising in random matrix theory. 

One special property of this distribution function is that 
it has mean zero. 

\begin{prop}\label{prop-mean1}
  We have 
\begin{equation} 
 \int_{-\infty}^\infty xd\F0(x)=0.
\end{equation} 
\end{prop}

\begin{proof}
We note that the term in front of $v(x)$ in the definition 
of $\F0$ has another expression : 
\begin{equation} 
  x+ 2u'(x)+2u(x)^2 = E(x)^{-4} \int_{-\infty}^x E(t)^4dt.
\end{equation} 
This follows by noting that 
$y(x):=x+2u'+2u^2$ satisfies 
\begin{equation} 
y'(x)=1+2u(x)y(x), 
\qquad 
y(x)=\frac1{\sqrt{-2x}}(1+o(1)), \quad x\to -\infty.
\end{equation} 

Then we have 
\begin{equation} 
  \F0(x) = {d\over dx} 
\biggl\{ \FGUE(x) \int_{-\infty}^x E(t)^4 dt \biggl\}, 
\end{equation} 
which upon integrating gives
\begin{equation}\label{en} 
\int_{-\infty}^x \F0(t) dt
=
\FGUE(x) \int_{-\infty}^x E(t)^4 dt
=
\FGUE(x) E(x)^4 \bigl( x + 2 u'(x) + 2 u(x)^2 \bigr).
\end{equation} 
Since 
$\int_{-\infty}^\infty xd\F0(x)=
\lim_{x\to\infty} \bigl[ x\F0(x)-\int_{-\infty}^x \F0(y)dy \bigr]$ 
from integration by parts,
subtracting \eqref{en} from $x \F0(x)$ and taking the limit $x\to\infty$, 
we find that $\F0$ has mean $0$, as required.
\end{proof}

\begin{rem}
  The mean zero property of $\F0(x)$ was suggested in \cite{SpohnP1} 
by numerical computation. 
Moreover, by an indirect argument for PNG with $\alpha_\pm=1$,
it is shown that the average of $L(t)$ is
exactly $2t$ (see \cite{SpohnP3}), 
which implies that $\F0$ has mean zero.
We note that the means of $\FGUE$ and $\FGOE$
are $-1.77109\cdots$ and $-0.76007\cdots$, respectively.
\end{rem}

\section{Around the transition : $\alpha_\pm\to 1$}\label{sec:scale}

In this section, we investigate the 
transition around $\alpha_\pm=1$ in detail.


To state the results, Theorem \ref{thm14}, we need some preliminary 
definitions.
Let $\Gamma$ be the real line $\R$, oriented from $+\infty$ to $-\infty$.
Let $m(\cdot;x)$ be the solution to the Painlev\'e II Riemann-Hilbert 
problem : 
\begin{equation}\label{as20}
 \begin{cases}
    m(z;x) \qquad \text{is analytic in $z\in\C\setminus\Gamma$,}\\
    m_+(z;x)=m_-(z;x) \begin{pmatrix} 1& -e^{-2i(\frac43z^3+xz)}\\
e^{2i(\frac43z^3+xz)}& 0 \end{pmatrix} \quad \text{for $z\in\Gamma$,}\\
    m(z;x) = I+O\bigl(\frac1{z}\bigr) \qquad \text{as $z\to\infty$.}
 \end{cases}
\end{equation}
Here $m_+(z;x)$ (resp., $m_-$) is the limit of $m(z';x)$ as $z'\to z$
from the left (resp., right) of the contour $\Gamma$ :
$m_\pm(z;x)=\lim_{\epsilon\downarrow 0}m(z\mp i\epsilon;x)$.
The relation between the above Riemann-Hilbert problem 
and the Painlev\'e II equation is the following 
(see, e.g., \cite{DZ2}).
If we expand 
\begin{equation}\label{as7.27}
   m(z;x) = I+ \frac{m_1(x)}{z} + O\bigl(\frac1{z^2}\bigr),
\qquad \text{as $z\to\infty$},
\end{equation}
we have
\begin{eqnarray}
\label{as22}   2i(m_1(x))_{12} = -2i(m_1(x))_{21} &=& u(x), \\
\label{as23}   2i(m_1(x))_{22} = -2i(m_1(x))_{11}&=& v(x),
\end{eqnarray}
where $u(x)$ and $v(x)$ are defined in \eqref{as10}-\eqref{as12}.  The
above Riemann-Hilbert problem is the special case of monodromy data
$p=-q=1$, $r=0$ in the standard family of Painlev\'e II Riemann-Hilbert
problems.

Define 
\begin{eqnarray}
  a(x,w)&=&\begin{cases}
m_{22}(-iw;x) &w>0,\\
-m_{21}(-iw;x)e^{\frac83w^3-2xw} &w<0,
\end{cases}\\
  b(x,w)&=&\begin{cases}
m_{12}(-iw;x) &w>0,\\
-m_{11}(-iw;x)e^{\frac83w^3-2xw} &w<0.
\end{cases}
\end{eqnarray}
From the jump condition of the Riemann-Hilbert problem \eqref{as20}, 
$a(x,w)$ and $b(x,w)$ are continuous at $w=0$.
Indeed since $m_-v$ in the upper half plane of $\C$ is an analytic 
continuation of $m_+$ in the lower half plane of $\C$, 
$a(x,w)$ and $b(x,w)$ are analytic in $w$.
We have the following properties of $a,b$.

\begin{lem}\label{lem13}
For all $x,w\in\R$, we have :
\begin{enumerate}
\item $a(x,w), b(x,w)$ are real.
\item For each fixed $w\in\R$, 
\begin{eqnarray}
  a(x,w)&=& I+O(e^{-cx^{3/2}}), \qquad\qquad x\to +\infty,\\
  b(x,w)&=& -e^{\frac83w^3-2xw} \bigl( I+O(e^{-cx^{3/2}})\bigr), 
\qquad x\to+\infty,\\
  a(x,w) &\sim& 
\frac1{\sqrt{2}} e^{\frac43 w^3-\frac{\sqrt{2}}{3}|x|^{3/2}+|x|w
-\sqrt{2}w^2|x|^{1/2}}
\qquad x\to-\infty,\\
  b(x,w) &\sim&
-\frac1{\sqrt{2}} e^{\frac43 w^3-\frac{\sqrt{2}}{3}|x|^{3/2}+|x|w
-\sqrt{2}w^2|x|^{1/2}}
\qquad x\to-\infty.
\end{eqnarray}
\item
\begin{eqnarray}
\lim_{w\to+\infty} a(x,w)=1, && \lim_{w\to+\infty} b(x,w)=0,\\
\lim_{w\to-\infty} a(x,w)=0, && \lim_{w\to-\infty} b(x,w)=0,\\
a(x,0)=E(x)^2, && b(x,0)=-E(x)^2.
\end{eqnarray}
\item
\begin{eqnarray}
  \frac{\partial}{\partial x} a(x,w) &=& u(x)b(x,w), \\
\frac{\partial}{\partial x} b(x,w) &=& u(x)a(x,w) -2w b(x,w),\\
  \frac{\partial}{\partial w} a(x,w) &=& 
2(u(x))^2a(x,w)-\bigl(4wu(x)+2u'(x)\bigr)b(x,w), \\
\frac{\partial}{\partial w} b(x,w) &=& 
\bigl(-4wu(x)+2u'(x)\bigr) a(x,w)+\bigl(8w^2-2x-2(u(x))^2\bigr)b(x,w).
\end{eqnarray}
\item
\begin{eqnarray} 
  a(x,w)&=&-b(x,-w) e^{\frac83w^3-2xw}, \\
  b(x,w)&=&-a(x,-w) e^{\frac83w^3-2xw}.
\end{eqnarray} 
\item For each fixed $y\in\R$, as $w\to-\infty$, 
\begin{eqnarray} 
  a(2y\sqrt{|w|}+4w^2,w) &\to& \erf(y),\\
b(2y\sqrt{|w|}+4w^2,w) &\sim& -e^{\frac{16}3|w|^3+4y|w|^{3/2}}.
\end{eqnarray} 
\end{enumerate}
\end{lem}

\begin{proof}
The properties (i)-(v) are consequences of
\eqref{as20} and Lemma 2.1 in \cite{BR2}.
The result (vi) is obtained by applying the Deift-Zhou method 
to the Riemann-Hilbert problem \eqref{as20}.
The specialty of the scaling $x=2y\sqrt{|w|}+4w^2$ is 
related to the fact that 
the exponent term $\frac43z^3+xz$ of the anti-diagonal entry 
in the jump matrix of the Riemann-Hilbert problem \eqref{as20} 
has the critical points at $z=\pm i\frac{\sqrt{x}}{2}$, which 
are the stationary phase points in the asymptotic analysis 
(see \cite{DZ2}).
Recalling $z=-iw$, we see that $x=4w^2$ corresponds to one of these 
stationary phase points.
By analyzing $m(-iw;x)$ around this stationary phase point, 
we obtain $m_{12}(-iw;x)e^{\frac83w^3-2xw}\sim -\erf(y)$ 
and $m_{11}(-iw;x) \sim 1$ as $w\to -\infty$ with the above $x$.
Similar computation appeared in Section 10.3 of \cite{BR2}, 
and we omit the detailed computations here.
\end{proof}

We now define the following functions with parameters $w_+, w_-, w$.

\begin{defn}
For each $w_+,w_-\in\R$, when $w_++w_-\neq 0$, set 
\begin{eqnarray}
  H(x;w_+,w_-)&=&\biggl\{ a(x,w_+)a(x,w_-)   
-\frac{a(x,w_+)a(x,w_-)-b(x,w_+)b(x,w_-)}{2(w_++w_-)}
v(x)\biggr\} \FGUE(x).
\end{eqnarray}
When $w_++w_-=0$, we use the l'Hopital's 
rule (note Lemma \ref{lem13} (v)).
Also using Lemma \ref{lem13} (iii), set 
\begin{equation}
  G(x;w)=\lim_{w_-\to+\infty} H(x;w,w_-)
= a(x,w) \FGUE(x).
\end{equation}
\end{defn}

From Lemma \ref{lem13} (ii), $H$ and $G$ have the limit $1$ as $x\to\infty$ 
and has the limit $0$ as $x\to-\infty$ for each fixed $w_+,w_-, w$.
Also theorem \ref{thm14} below shows that they are limits of distribution 
functions. 
Therefore $H$ and $G$ are distribution functions.
These distribution functions interpolate between the functions
$\FGUE$, $\FGOE^2$ and $\F0$ of theorem \ref{thm1}.
The following results follow from Lemma \ref{lem13} (iii), (vi).

\begin{prop}\label{prop13}
For fixed $x\in\R$, we have 
\begin{equation}\label{a1}
  H(x;w_+,w_-)\to \begin{cases}
\FGUE(x), \quad &w_+,w_-\to +\infty\\
\FGOE(x)^2, \quad &\text{$w_+=0, w_-\to +\infty$ or 
$w_-=0, w_+\to +\infty$}\\
\F0, \quad &w_+=w_-=0\\
0, \quad &\text{$w_+$ or $w_-$ $\to -\infty$},
\end{cases}
\end{equation}
and 
\begin{equation}\label{a2}
  G(x;w)\to \begin{cases}
\FGUE(x), \quad &w\to +\infty\\
\FGOE(x)^2, \quad &w=0\\
0, \quad &w\to -\infty.
\end{cases}
\end{equation}
Also we have for fixed $x\in\R$, with $w=-\max\{-w_+,-w_-\}$, 
as $w\to-\infty$, 
\begin{eqnarray} \label{a3}
   H(2x\sqrt{|w|}+4w^2;w_+,w_-) &\to& 
\begin{cases} \erf(x), \quad &w_+\neq w_-\\
\erf^2(x), &w_+=w_-, 
\end{cases}\\
 \label{a4}  
G(2x\sqrt{|w|}+4w^2;w) &\to& \erf(x). 
\end{eqnarray} 
\end{prop}

Now the main theorem in this section is that if we take proper 
scaling of $\alpha_\pm\to 1$ in $t$, we obtain the above 
functions in the limit.

\begin{thm}\label{thm14}
Set $w_\pm$ by 
\begin{equation}
  \alpha_\pm=1-\frac{2w_\pm}{t^{1/3}}.
\end{equation}
We have the followings.
\begin{enumerate}
\item 
  When $0\le \alpha_+<1$ and $w_-\in\R$ are fixed, 
\begin{equation}   
  \lim_{t\to\infty}  \Prob\biggr( \frac{L(t)-2t}{t^{1/3}}\le x \biggr) =
G(x;w_-).
\end{equation} 
When $w_+\in\R$ and $0\le \alpha_-<1$ are fixed, the limit 
is $G(x;w_+)$.
\item 
When $w_\pm\in\R$ are fixed, as $t\to\infty$,
\begin{equation}   
  \lim_{t\to\infty} \Prob\biggr( \frac{L(t)-2t}{t^{1/3}}\le x \biggr) =
H(x;w_+,w_-).
\end{equation} 
\end{enumerate}
\end{thm}

A special case of (ii) is when $w_+=-w_-$, which corresponds to 
$\alpha_+\alpha_-=1$.
In this case, the limiting shape of PNG has curvature $0$ 
(see \cite{SpohnP1}).
Hence we obtain a one-parameter family of distribution
functions for fluctuations of a flat curvature PNG.
We have the exact values of the means which include 
Proposition \ref{prop-mean1} as a special case when $w=0$.

\begin{prop}\label{prop-mean2}
We have for each $w\in\R$, 
\begin{equation} 
  \int_{-\infty}^\infty xdH(x;w,-w)=4w^2.
\end{equation} 
\end{prop}

\begin{proof}
By l'Hopital's rule and Lemma \ref{lem13} (iv), (v), we have 
\begin{equation} 
  H(x;w,-w)= \bigl\{ a(x;w)a(x;-w) - y(x)v(x)\bigr\} \FGUE(x),
\end{equation} 
where 
\begin{equation} 
   y(x)= (2u^2+x-4w^2)a(x;w)a(x;-w)
-(u'+2wu)b(x;w)a(x;-w)
-(u'-2wu)a(x;w)b(x;-w).
\end{equation} 
Then by Lemma \ref{lem13} (iv), (v) again, we obtain 
\begin{equation}
   y'(x)= a(x;w)a(x;-w), 
\end{equation}
which implies that 
\begin{equation} 
  H(x;w,-w)= \bigl\{ y'(x) - y(x)v(x)\bigr\} \FGUE(x)
=\frac{d}{dx} \biggl\{ y(x)\FGUE(x) \biggr\}.
\end{equation} 
Hence as in the proof of Proposition \ref{prop-mean1}, 
\begin{equation} 
  \lim_{x\to\infty} 
\biggl[ xH(x;w,-w)-\int_{-\infty}^x H(y;w,-w)dy \biggr]
= \lim_{x\to\infty} (xy'-xvy-y)\FGUE(x) = 4w^2, 
\end{equation} 
using the asymptotics of $u, v, a, b$.
\end{proof}


\begin{rem}\label{rem1}
The function $G(x;w)$ appeared 
as $F^\symmu(x;w)$ in (2.21), (2.22) of \cite{BR2}. 
Indeed this function was obtained 
in a different point selection process.
Namely, in the open square $(0,1)\times (0,1)$, 
let $\delta=\{(t,t):0<t<1\}$, the diagonal line, 
and let $\delta^t=\{(t,1-t): 0<t<1\}$, the anti-diagonal line.
We select $4P(t^2)$ points in $(0,1)\times (0,1) \setminus 
(\delta\cup\delta^t)$, 
$2P(\alpha t)$ points on $\delta\setminus(\frac12,\frac12)$,
and $2P(\beta t)$ points on $\delta^t\setminus(\frac12,\frac12)$
such that the resulting point configuration is symmetric with respect to
both $\delta$ and $\delta^t$.
Let $L^\symmu(t;\alpha,\beta)$ be the length of the longest 
up/right path of a random point configuration in this point selection process.
In \cite{BR2}, it is proved that 
for \emph{any} fixed $\beta\ge 0$, 
\begin{eqnarray}
  \Prob\biggr( \frac{L^\symmu(t;\alpha,\beta)-4t}{2t^{1/3}}\le x \biggr) &\to&
\begin{cases}
  \FGUE(x), \quad &0\le \alpha <1,\\
  \FGOE(x)^2, \quad &\text{$\alpha=1$,}
\end{cases}\\
  \Prob\biggr( \frac{L^\symmu(t;\alpha,\beta)-2(\alpha+\alpha^{-1})t}
{\sqrt{2(\alpha-\alpha^{-1})}t^{1/2}}\le x \biggr) &\to&
\erf(x), \qquad \alpha>1,\\
  \Prob\biggr( \frac{L^\symmu(t;\alpha,\beta)-4t}{2t^{1/3}}\le x \biggr) &\to&
G(x;w), \qquad \alpha=1-\frac{2w}{t^{1/3}}
\end{eqnarray}
This result is identical to the limits of $L(t)$ when $\alpha_-<1$ 
is fixed.
Comparing this with theorems \ref{thm1} and \ref{thm14}, 
we see that after scaling, 
$L^\symmu(t;\alpha,\beta)$ with any fixed $\beta$
and $L(t)$ with any fixed $\alpha_-<1$ have the same statistics in 
the limit $t\to\infty$.
Indeed, one can prove more than that.
In finite $t$, $L(2t)$ with $\alpha_-=0$ and $L^\symmu(t;\alpha_+,0)$ 
are the same (see Remark at the end of Step 1. in Section \ref{sec:proof}.)
\end{rem}

\section{Lattice directed polymer and exclusion process}\label{sec:excl}

In this section, we analyze a certain lattice directed polymer problem 
which is closely related to the point selection model discussed above.

Let $\N=\{1,2,3,\cdots\}$ and $\N^*=\{0,1,2,\cdots\}$.
Let $0<q<1$, $\alpha_\pm\ge 0$ be fixed numbers 
such that $\alpha_\pm\sqrt{q}<1$.
We denote by $g(q)$ the geometric distribution with 
parameter $q$ : 
\begin{equation}  
   \Prob(g(q)=k)=(1-q)q^k, \qquad k=0,1,2,\cdots.
\end{equation} 
When a random variable $w$ has distribution $g(q)$, 
we use the notation $w\sim g(q)$.
At each site $(i,j)\in \N^*\times \N^*$, 
we attach a random variable $w(i,j)$ where 
\begin{eqnarray} 
  w(i,j) &\sim& g(q), \qquad (i,j)\in \N\times \N,\\
  w(i,0) &\sim& g(\alpha_+\sqrt{q}), \qquad i\in \N,\\
  w(0,j) &\sim& g(\alpha_-\sqrt{q}), \qquad j\in \N,\\
  w(0,0) &=& 0. 
\end{eqnarray}  
We call a collection $\pi$ of sites in $\N^*\times \N^*$ 
an up/right path 
if when $(i,j)\in\pi$, either 
$(i+1,j)\in\pi$, or $(i,j+1)\in\pi$.
Let $Path(N)$ be the set of all up/right paths from $(0,0)$ 
to $(N,N)$.
Define 
\begin{equation}\label{eq-excl6}
  X(N) = \max \{ \sum_{(i,j)\in \pi} w(i,j) : \pi\in Path(N) \}.
\end{equation} 

The special case $\alpha_+=\alpha_-=0$ was introduced by Johansson 
in \cite{kurtj:shape}; the above model adds a special row 
and column to his model. 
In \cite{kurtj:shape}, it is proved that 
\begin{equation}  
  \lim_{N\to\infty}\Prob\biggl(\frac{X(N)-\mu(q)N}{\sigma(q)N^{1/3}}
\le x \biggr)=\FGUE(x), \qquad \alpha_+=\alpha_-=0,
\end{equation} 
where 
\begin{equation}\label{e-excl8}  
  \mu(q)=\frac{2\sqrt{q}}{1-\sqrt{q}}, 
\qquad \sigma(q)=\frac{q^{1/6}(1+\sqrt{q})^{1/3}}{1-\sqrt{q}}.
\end{equation} 

It is shown in \cite{kurtj:shape} that 
this directed polymer model can be interpreted 
as a growth model in 2-dimensional 
space, or as a discrete exclusion process.
The corresponding exclusion process in our case is the following.
We use the notation $+$ for the location of a particle.
If a site is vacant, we use the notation $-$.
Initially there are particles at the sites 
$\{\cdots, -4, -3,-2\}\cup\{0\}\subset \Z$. 
Hence the initial configuration on $\Z$ can be written as 
$(\cdots, +,+,+,-,+,-,-,-,\cdots)$ where the leftmost $-$ is at the 
site $-1$ and the rightmost $+$ is at the site $0$.
At each (discrete) time, the rightmost particle jumps to its right site 
with probability $1-\alpha_+\sqrt{q}$, 
and the leftmost \emph{hole} jumps to its left with probability 
$1-\alpha_-\sqrt{q}$, while in the `bulk', each particle 
jumps to its right (equivalently, a hole jumps to its left) 
with probability $1-q$ if its right site is vacant.

\medskip
As in the point selection model, when $\alpha_\pm$ are small enough, 
$X(N)$ would have the limiting distribution $\FGUE$ 
as in the case when $\alpha_\pm=0$. 
For general $\alpha_\pm$, we obtain
results parallel to those in Theorems \ref{thm1} and \ref{thm14}.
After the following changes, we obtain the same limiting results 
as in Theorems \ref{thm1} and \ref{thm14}. 
\begin{enumerate}
\item Every limit is taken as $N\to\infty$.
\item The scaled random variable in Theorem \ref{thm1} (i) and 
Theorem \ref{thm14} 
is now 
\begin{equation}\label{enew1}  
  \frac{X(N)-\mu(q)N}{\sigma(q)N^{1/3}}
\end{equation} 
where $\mu(q)$ and $\sigma(q)$ are defined in \eqref{e-excl8}.
\item
The scaling of $\alpha_\pm$ in Theorem \ref{thm14} is now 
\begin{equation}\label{enew2} 
  \alpha_\pm= 1 - \frac{2w_\pm}{\sigma(q)N^{1/3}}.
\end{equation} 
\item  
The scaling in Theorem \ref{thm1} (ii) is now 
\begin{equation}\label{enew3} 
  \frac{X(N)-\eta(q,\alpha)N}{\rho(q,\alpha)\sqrt{N}}, 
\qquad 1<\alpha <\frac1{\sqrt{q}},
\end{equation} 
where
\begin{equation} 
\eta(q,\alpha)=\frac{\sqrt{q}(\alpha+\alpha^{-1}-2\sqrt{q})}
{(1-\sqrt{q}\alpha)(1-\sqrt{q}\alpha^{-1})}, 
\qquad 
\rho(q,\alpha)=\frac{\sqrt{q}(\alpha-\alpha^{-1})^{1/2}
(\sqrt{q}^{-1}-\sqrt{q})^{1/2}}
{(1-\sqrt{q}\alpha)(1-\sqrt{q}\alpha^{-1})}.
\end{equation} 
We note that $\eta(q,\alpha)> \mu(q)$ for $\alpha>1$.
\end{enumerate}


\begin{rems}
In the above model (also in \cite{BR2}), 
only the paths ending at the diagonal point $(N,N)$ are considered.
It is of interest to obtain similar results for path ending at 
general point $(M,N)$.
The difficulty in this general case comes from the fact that 
we need asymptotics of orthogonal polynomials with 
respect to the weight function $(1+\sqrt{q}z)^M(1+\sqrt{q}z^{-1})^N$ 
(see Step 2. of Section \ref{sec:proof}).
This weight function is not real for $|z|=1$, 
which makes the Riemann-Hilbert method (which we employed to obtain 
asymptotics) more difficult to analyze.
We are planning to come back to this problem in the future.

We note that in \cite{kurtj:shape}, Johansson was able to 
obtain results for this general case when $\alpha_+=\alpha_-=0$.
He have used different determinant expression (of 
Fredholm type rather than Toeplitz type) involving 
different orthogonal polynomials (which is discrete), 
and did not involve the non-real weight function.
\end{rems}

\begin{rems}
  If we take the limit as $q\to 1$, we obtain exponential random variables 
instead of geometric distribution 
(see \cite{kurtj:shape}). 
In order to investigate the exponential random variables 
case, we set $q=1-\frac1{L}$ and $l=xL$, and let $L\to\infty$ 
in the determinantal formula \eqref{e-excl9}-\eqref{e-excl10}.
This double scaling limit is not carried out yet, 
which we plan to do in the future.
\end{rems}

\begin{rems}
  If  we take the limit as $q\to 0$, we have rare events, hence we obtain 
the Poisson process discussed above.
Indeed, by setting $\sqrt{q}=\frac{t}{N}$, and letting $N\to\infty$, 
we recover the Theorems \ref{thm1} and \ref{thm14} 
from \eqref{enew1}-\eqref{enew3}.
\end{rems}

\section{Proofs}\label{sec:proof}

In this section, we sketch how one can obtain 
Theorems \ref{thm1} and \ref{thm14} using the results of 
\cite{BR1,BR2}.
We split the proof into three steps.

\bigskip
{\bf Step 1.}
We will prove the following formula : 
\begin{equation}\label{e-f1}
  \Prob(L(t)\le l) =
e^{-(\alpha_++\alpha_-)t-t^2}
(D_l' - \alpha_+\alpha_-D'_{l-1})
\end{equation}
where 
\begin{equation}\label{e-f2} 
  D'_l=
\Exp_{U\in U(l)} \det  
(1+\alpha_+U)(1+\alpha_-U^\dagger)
e^{2t \Re\Tr(U) }.
\end{equation} 

Suppose $\alpha_+\alpha_-<1$.
In the lattice directed polymer model discussed in Section \ref{sec:excl}, 
add a random variable $w(0,0)\sim g(\alpha_+\alpha_-)$ 
at the site $(0,0)$. 
Let $X^+(N)$ be defined by the formula \eqref{eq-excl6} 
with this new random variable added. 
Now this is a special case of the model $\tsymmU$ of 
(7.5)-(7.7) discussed in 
Section 7 of \cite{BR1} : 
take $W=W'=\{0,1,\cdots,N\}\subset \N^*$, 
and take $q_0=\alpha_-$, $q_0'=\alpha_+$ and 
$q_j=q_j'=\sqrt{q}$ for $1\le j\le N$ 
(in \cite{BR1}, we have taken $W,W'\subset \N$, but 
simple translation makes no change.)
From Theorem 7.1 (7.30) of \cite{BR1}, we obtain 
\begin{equation}\label{eq-pf1}
\begin{split}
  \Prob(X^+(N)\le l) = &(1-\alpha_+\alpha_-)(1-\alpha_+\sqrt{q})^N
(1-\alpha_-\sqrt{q})^N(1-q)^{N^2}\\
&\cdot \Exp_{U\in U(l)} \det\bigl\{
(1+\alpha_+U)(1+\alpha_-U^\dagger)
(1+\sqrt{q}U)^N(1+\sqrt{q}U^\dagger)^N
\bigr\}.
\end{split}
\end{equation} 
Now we set $\sqrt{q}=\frac{t}{N}$, and take the limit $N\to\infty$. 
Then inside the square $\{1,\cdots,N\}^2$, 
we obtain a Poisson process of parameter $t^2$ : 
$\Prob(w=0)=1-\frac{t^2}{N^2}$, $\Prob(w=1)=(1-\frac{t^2}{N^2})\frac{t^2}{N^2}$
and $\Prob(w\ge 2)=\frac{t^4}{N^4}$.
The probability that there are $k$ points in the square
$\{1,\cdots,N\}^2$ is equal to 
\begin{equation}  
   \binom{N}{k}\biggl(1-\frac{t^2}{N^2}\biggr)^{N-k}
\biggl\{(1-\frac{t^2}{N^2})\frac{t^2}{N^2}\biggr\}^k 
\to \frac{e^{-t^2}(t^2)^k}{k!},
\end{equation} 
disregarding the events of having more than one points 
at one site which has probability $0$ in the limit.
Similarly we obtain Poisson process with parameter $\alpha_+t$ on the 
bottom edge (not including the origin), and Poisson process 
with parameter $\alpha_-t$ on the left edge (not including the origin).
At the origin, we have a geometric distribution with 
parameter $\alpha_+\alpha_-$.
Let $L^+(t)$ be the length of the longest up/right path 
in this process.
This is related to $L(t)$ by
\begin{equation}\label{eq-pf3}
  L^+(t)=L(t)+\chi, \qquad \chi\sim g(\alpha_+\alpha_-),
\end{equation}
since any up/right path will include the $g(\alpha_+\alpha_-)$ points in 
the lower-left corner.
From \eqref{eq-pf1}, we obtain 
\begin{equation}  
  \Prob(L^+(t)\le l) =
(1-\alpha_+\alpha_-) e^{-(\alpha_++\alpha_-)t-t^2} 
D_l'
\end{equation} 
where $D_l'$ is defined in \eqref{e-f1}.

In order to obtain the formula for $L(t)$, set
\begin{equation}  
  Q(x)=\sum_{l\ge 0} \Prob(L(t)\le l)x^l, 
\qquad Q^+(x)=\sum_{l\ge 0} \Prob(L^+(t)\le l)x^l.
\end{equation} 
Then using \eqref{eq-pf3}, we obtain 
\begin{equation}  
  Q^+(x)=\sum_{0\le k\le l} \Prob(L(t)\le l-k)\Prob(\chi=k)
= (1-\alpha_+\alpha_-)(1-\alpha_+\alpha_- x)^{-1} Q(x).
\end{equation} 
Thus by comparing the coefficients, we have \eqref{e-f1} 
for $\alpha_+\alpha<1$.

Observe that the first quantity in the right-hand side in \eqref{e-f1} 
is entire in $\alpha_+, \alpha_-$, and the second quantity is 
polynomial in $\alpha_+, \alpha_-$, hence is entire.
Since both sides of \eqref{e-f1} agree analytically for $\alpha_+\alpha_-<1$ 
and the right-hand side is entire, 
they agree in general $\alpha_\pm$ where they both converge and are defined.
Thus \eqref{e-f1} holds 
for $0\le \alpha_\pm$.
Similar consideration yields the following formulae for the probability of 
the lattice directed polymer problem :
\begin{equation}\label{e-excl9}
 \Prob(X(N)\le l)
=
(1-\alpha_+\sqrt{q})^N(1-\alpha_-\sqrt{q})^N(1-q)^{N^2}
(T'_l - \alpha_+ \alpha_- T'_{l-1})
\end{equation}
where
\begin{equation}\label{e-excl10}
T'_l = E_{U\in U(l)} \det \{ (1+\alpha_+ U)(1+\alpha_-U^\dagger)
(1+\sqrt{q} U)^N(1+\sqrt{q}U^\dagger)^N\},
\end{equation}
and $\alpha_\pm$ are subject to the constraint
$0\le \alpha_{\pm}<1/\sqrt{q}$.

\begin{rem}\label{rem2}
It is the formula \eqref{e-f1} that makes connection 
with the process $\tsymmu$ 
mentioned in Remark of Section \ref{sec:scale}.
In \cite{BR1} (4.16), it is proved that 
$\Prob(L^\symmu(t;\alpha,0)\le 2l)$ is given by \eqref{e-f1} 
with $\alpha_+=\alpha$ and $\alpha_-=0$.
It is not clear why these two processes should be the same.
\end{rem}

\bigskip
{\bf Step 2.}
Let
\begin{equation}\label{e-int6}
D_l = E_{U\in U(l)} e^{2t\Re\Tr(U)},
\end{equation}
which is in another form, the $l\times l$ Toeplitz determinant
$\det(c_{j-k})_{0\le j,j <l}$ where $c_j$ is the Fourier coefficient of
$e^{2t\cos\theta}$.  In \cite{BDJ}, the asymptotics of $D_l$ in the double
scaling limit as $t,l\to\infty$ was studied; hence the second step in the
proof of Theorem \ref{thm1} is to eliminate the term $\det(1+\alpha_+
U)(1+\alpha_-U^\dagger)$ from the integrand in \eqref{e-f1}.
This step is established in Theorem 3.2 of \cite{BR1}.
Let $\pi_n(z)=z^n+\cdots$ be the $n^{\text{th}}$ monic orthogonal polynomial
with respect to the weight $e^{t(z+z^{-1})} dz/(2\pi iz)$ on the
unit circle $|z|=1$ where $z\in\C$ :
\begin{equation}\label{e-in8}
   \int_{|z|=1} \pi_n(z)\overline{\pi_m(z)} e^{t(z+z^{-1})} \frac{dz}{2\pi iz}
= \delta_{nm} N_n,
\end{equation}
for some constants $N_n$.
Note that since the weight function is real, all coefficients
of $\pi_n$ are real.
We define
\begin{equation}\label{e-in9}
   \pi^*_n(z)=z^n\pi(z^{-1}).
\end{equation}
In Theorem 3.2 of \cite{BR1},
using the Weyl integration formula for $U(l)$ and
familiar Vandermonde type argument
together with the relations between orthogonal polynomials
on the unit circle and those on the interval,
it is proved that
for $\alpha_+\alpha_-\neq 1$,
\begin{equation}\label{e-in10}
D'_l = 
\frac{ \pi^*_l(-\alpha_+) \pi^*_l(-\alpha_-)
- \alpha_+\alpha_- \pi_l(-\alpha_+)\pi_l(-\alpha_-) }  
{1-\alpha_+\alpha_-}
D_l 
\end{equation}
For $\alpha_+\alpha_-=1$, l'Hopital's rule applies implying
with $\alpha_+=\alpha$, $\alpha_-=1/\alpha$,
\begin{equation}\label{e-in11}
D'_l = \bigl\{ (1-l)\pi_l(-\alpha)\pi_l(-\alpha^{-1})
- a \pi'_l(-\alpha) \pi_l(-\alpha^{-1})
- \alpha^{-1}\pi_l(-\alpha) \pi'_l(-\alpha^{-1}) \bigr\} D_l.
\end{equation}
We obtain similar formulae for the lattice directed polymer problem: 
replace the weight $e^{t(z+z^{-1})}$ by 
$(1+\sqrt{q}z)^N(1+\sqrt{q}z^{-1})^N$.

\bigskip
{\bf Step 3.}
The remaining task is to obtain the asymptotics of $D_l$ and $\pi_l(-\alpha)$ 
as $l,N\to\infty$ in a proper rate.
This is obtained in \cite{BDJ, BR2} by applying the Deift-Zhou steepest 
descent method (see \cite{DZ1, deift}) to the Riemann-Hilbert problem 
for the orthogonal polynomials $\pi_l(z)$.
These asymptotic results are summarized in Section 5 of \cite{BR2}. 
Theorem \ref{thm1} follows by plugging in these asymptotics into 
\eqref{e-in10}, \eqref{e-in11};
we omit the calculations.
There are similar asymptotic results for the lattice directed polymer problem 
which yield Theorem \ref{thm14};
see Proposition 3.2 and subsequent remarks in \cite{B00}.

\bibliographystyle{plain}
\bibliography{paper6}

\end{document}